\documentclass[12pt]{amsart}
\usepackage{amsmath, amssymb}
\usepackage{graphicx}
\usepackage{tikz}
\usetikzlibrary{backgrounds,patterns,shapes,calc}

\newtheorem{theorem}{Theorem}
\newtheorem{conjecture}[theorem]{Conjecture}
\newtheorem{observation}[theorem]{Observation}
\newtheorem{corollary}[theorem]{Corollary}
\newtheorem{proposition}[theorem]{Proposition}
\newtheorem{lemma}[theorem]{Lemma}

\newtheorem{definition}[theorem]{Definition}

\newtheorem{quest}{Question}

\def\bth{\begin{theorem}}
\def\eth{\end{theorem}}
\def\bc{\begin{corollary}}
\def\ec{\end{corollary}}
\def\bcj{\begin{conjecture}}
\def\ecj{\end{conjecture}}

\def\co{{\rm co}}

\newcommand{\bizveg}{{\hfill $\Box$}}

\newcommand{\ie}{{ i.e.\ }}

\newcommand{\whp}{{\bf whp}}

\def\Pr{\mathbb{P}}

\title[Clustering with coloring]{Graph Clustering via Generalized Colorings}

\author[London]{Andr\'{a}s London}\address{Department of Computer Science,
University of Szeged and Department of Operations Research, Pozna\`{n} University of Economics}

\thanks{The first author's research is partially supported by the National Research, Development and Innovation Office – NKFIH, Fund No. SNN-135643.}
\email{london@inf.u-szeged.hu}

\author[Martin]{Ryan R. Martin}\address{Department of Mathematics, Iowa State University}

\thanks{The second author's research is partially supported by Simons Foundation Collaboration Grant \#353292 and by the J. William Fulbright Educational Exchange Program.}
\email{rymartin@iastate.edu}

\author[Pluh\'ar]{Andr\'as Pluh\'ar}\address{Department of Computer Science,
University of Szeged}

\thanks{The third author's research is partially supported by the project “Deepening the activities of the
Hungarian Industrial Innovation Mathematical Service Network HUMATHS-IN,” No. EFOP-3.6.2-16-2017-00015.
}
\email{pluhar@inf.u-szeged.hu}


\begin{document}

\begin{abstract}
We propose a new approach for defining and searching clusters in graphs that represent real technological or transaction networks.
In contrast to the standard way of finding dense parts of a graph, we concentrate on the structure of edges 
between the clusters, as it is motivated by some earlier observations, e.g. in the structure of networks in ecology and economics and by applications of discrete tomography. Mathematically special colorings and chromatic numbers of graphs are studied.      
\end{abstract}

\maketitle

\section{Introduction and Results}

One of the main tasks in network theory is clustering vertices, see Newman
\cite{NG}. Graph clustering is a well-studied problem and has important applications in 
graph mining or model construction. The usual methods try to achieve many edges inside 
clusters and only a few between distinct clusters \cite{Sch}. This approach generally works well for so-called 
{\em social graphs}, which usually contain more triangles than a random graph with similar edge density or degree properties. 
In contrast {\em technological} or {\em transaction graphs} contain fewer triangles and often display tree-like structures. To measure the algorithms' efficacy the parameter known as {\em Newman modularity} is commonly used \cite{NG}.  

However, this standard approach is not always justified. Certain bipartite graphs, e.g. those that describe pollinator 
networks or trade networks, suggest the presence different structures, like the notion of 
{\em embeddedness}, see Uzzi \cite{Uzzi}. That is, the vertices of each color class can be ordered, and the smaller ranked vertex 
neighborhood contains the neighborhood of any higher ranked one. 
In the context of image processing, Junttila and Kaski \cite{Jun} call a binary matrix $A$ 
(that is, a matrix whose entries are either zero or one) {\em fully nested} if its rows and columns can be reordered such that the ones are in an echelon form. Let $G_A$ be the bipartite graph whose adjacency matrix is $A$. Then
$A$ being fully nested is equivalent to $G_A$ satisfying embeddedness. 

Let $X$ (the columns) and $Y$ (the rows) be the bipartition of $G_A$. 
The matrix $A$ and the graph $G_A$ are each said to be {\em $k$-nested with respect to $X$} if $X$ can be partitioned as $X_1, \dots, X_k$ such that all subgraphs spanned by $(X_i, Y)$ are fully nested for $i=1, \dots, k$.
The quantity of interest for any $G_A$ is smallest $k$ for which $G_A$ is $k$-nested. 
  
We present a new kind of clustering of general (that is, not necessary bipartite) transaction graphs 
via a certain class of proper colorings. The clusters are the color classes, since we do not want edges inside a cluster,
and we restrict the structure of the edges between the pairs of classes. The above examples 
suggest that in some cases there should be a fully nested or, equivalently, embeddedness relation among 
any two color classes. We generalize this notion to an arbitrary host graph $G$ and a forbidden
bipartite subgraph $H$ as follows.

\begin{definition}
Fix a bipartite graph $H$. A proper coloring of a graph $G$ is an \emph{$H$-avoiding coloring} if the union of any two color classes spans an induced $H$-free graph. 
Let $\chi_H(G)$ be the minimum number of colors in an $H$-avoiding coloring of $G$. 
\end{definition}

Note that the function $\chi_H(G)$ is not necessarily monotone either in $H$ or in $G$. However, 
we have a useful property: 

\begin{observation} \label{triv} For any graphs $H$ and $G$, $\chi(G) \leq\chi_H(G)$.
If $G$ is $H$-free, then $\chi(G)=\chi_H(G)$. 
\end{observation}

\subsection{Complexity issues}

We show that that the computation of $\chi_H(G)$ is NP-hard for some graphs, 
and polynomially computable for others. The most interesting case, when $H = 2K_2$,
gives back embeddedness as described above. For these generalized chromatic numbers we derive some
theoretical extremal results as well as results on complexity.\footnote{Heuristics for finding $H$-avoiding colorings and case studies will be presented in a future paper.}

In the following we use $K_n$,
$P_n$ and $C_n$ for the complete graph, path and cycle on $n$ vertices, respectively. For
graphs $H_1$ and $H_2$ on disjoint vertex sets, $H_1 \oplus H_2$ denotes their disjoint union.
Theorem~\ref{main_comp} gives a characterization of the complexity issues in computing
$\chi_H(G)$ depending on the graph $H$.

\begin{theorem}\label{main_comp} The computation of $\chi_H(G)$ is polynomial-time solvable 
if $H$ is $K_1 \oplus K_1$, $K_2$, or $K_2 \oplus K_1$ and is NP-hard for all other graphs.
\end{theorem}

It is valuable to spell out special cases since the proofs of these are needed in proving
Theorem~\ref{main_comp}.

\begin{lemma} \label{P_3} It is NP-complete to decide if $\chi_{P_3}(G) \leq 5$,
while it is polynomial time decidable if $\chi_{P_3}(G) \leq 3$.
\end{lemma}

\begin{lemma} \label{P_4} It is NP-complete to decide if $\chi_{P_4}(G) \leq 3$.
\end{lemma}

\begin{lemma} \label{K_2+K_1} There is a unique $H$-avoiding coloring of $G$ using
exactly $\chi_H(G)$ colors if $H=K_2 \oplus K_1$. One can find this coloring in
polynomial time.
\end{lemma}

Let us note that a $P_3$-avoiding coloring of $G$ has a nice combinatorial meaning,
it represents the edges of $G$ as the union of independent matchings. The computation
of $\chi_{P_3}(G)$ can be reduced to the normal chromatic number. Let $P_3(G)$ be a
graph which made from $G$ by adding an edge to every induced $P_3$, \ie making a
triangle out of these $P_3$.

\begin{observation} \label{P_3close}
$\chi(P_3(G))=\chi_{P_3}(G)$. 
\end{observation} 

Note that if a bipartite graph $G_A$ is $k$-nested then it has a similar reduction as 
in Observation~\ref{P_3close}.

For a bipartite graph $G_A$ with bipartition $(X,Y)$, define the {\em conflict graph} $\co(X)$ on $X$ such that
$(x, x')$ is an edge in $\co(G)$ for $x, x' \in X$ if there are $y, y' \in Y$ such that
$\{x, x', y, y'\}$ spans a $2K_2$ in $G_A$.   

\begin{observation} \label{k-nested}
The bipartite graph $G_A$ is exactly $k$-nested for $X$ if $\chi(\co(X))=k$. 
\end{observation}

For applications the computation of $\chi_{2K_2}(G)$ seems to be the most 
important case. 

\begin{theorem}\label{2K_2poly3} It is polynomial time decidable if $\chi_{2K_2}(G) \leq 3$. 
\end{theorem}

For a fix graph $H$ there is a linear upper bound on the value of $\chi_{H}(G)$. In this paper, $\lg n$ is the logarithm of $n$ in base 2 and $\log n$ is the natural logarithm of $n$.

\begin{proposition} \label{upper_bound}
	Let $H$ be a bipartite graph and let $k_1$ be the smallest positive integer such that each bipartition of $H$ has a part with size at least $k_1$. Let $k_2$ be the smallest positive integer such that each bipartition of $H$ has both parts of size at least $k_2$. Let $G$ be an $n$-vertex graph with chromatic number $\chi$ and independence number $\alpha$. If $k_2\geq 3$, then
	\begin{align*}
		\chi_{H}(G) 	\leq 	\min\left\{\frac{n}{k_1-1} + \frac{k_1-2}{k_1-1}\,\chi, \frac{n}{k_2-1}\left(1-\frac{1}{\chi}\right) + \frac{k_2-2}{k_2-1}(\chi-1) + 1\right\} .
	\end{align*}
	If $k_2=2$, then 
	\begin{align*}
		\chi_{H}(G) 	\leq 	\min\left\{\frac{n}{k_1-1} + \frac{k_1-2}{k_1-1}\,\chi, n-\alpha+1\right\} .
	\end{align*}
\end{proposition}

\subsection{Random graphs}

In the case where $G$ is a random graph drawn from $G(n,p)$, the Erd\H{o}s-R\'enyi random graph on $n$ vertices with edge probability $p$, we establish tight bounds for $\chi_{H}(G)$. The distribution of $\alpha(G)$, where $G\sim G(n,p)$ for $p$ fixed, was determined by Bollob\'as and Erd\H{o}s~\cite{BE}. The distribution of $\chi(G)$ was first proven by in a classic result by Bollob\'as~\cite{Bol} and the error terms have been further refined by various authors (see~\cite{AS}). 

To be precise, \whp~ means \emph{with high probability}, \ie a probability arbitrarily close to one, provided that the number of vertices (or other natural parameter) is large enough. 
\begin{theorem} \label{general_random}
	Let $H$ be a bipartite graph with $k_1$ and $k_2$ defined as in Proposition~\ref{upper_bound}. Fix $p\in (0,1)$, let $d=1/(1-p)$, and let $G\sim G(n,p)$. If $k_1\geq 3$ and $k_2\geq 2$, then there is a $C=C(H,p)$ such that \whp
	\begin{align*}
		\frac{n}{k_1-1} - C\log n 	\leq 	\chi_{H}(G) 	\leq 	\frac{n}{k_1-1} + O\left(\frac{n}{\log_d n}\right) .
	\end{align*}
	If $k_1=k_2=2$, then there exists a $C=C(H,p)$ such that \whp
	\begin{align*}
		n - C \log n	\leq 	\chi_{H}(G) 	\leq 	n - 2\log_d n + O\left(\log_d \log n\right) .
	\end{align*}
	In particular, if $H=2K_2$, then \whp
	\begin{align*}
		n - 8 \log_{1/Q} n + \Omega\left(\log_{1/Q} \log n\right) 	\leq 	\chi_{H}(G) 	\leq 	n - 2\log_d n + O\left(\log_d \log n\right) ,
	\end{align*}
	where $Q = 1 - 2 p^2 (1-p)^2$. 
\end{theorem}

Finally, we mention a useful observation on the $H$-avoiding chromatic number, see its
consequences in Section~\ref{corollaries}.

\begin{observation} \label{edge_bound}
	If $G$ is a graph such that every $H$-free induced subgraph has at most $\ell$ edges, 
	then $\chi_H(G)$ satisfies
	\begin{align*}
		\ell \binom{\chi_H(G)}{2} \geq e(G) .
	\end{align*}
\end{observation}

In the rest of the paper we provide the proofs of the main results. In Section~\ref{main} 
we prove Theorem~\ref{main_comp}. Section~\ref{polysec} contains the proof of Theorem~\ref{2K_2poly3}, 
while Section~\ref{random} contains the proofs of Proposition \ref{upper_bound} and of Theorem \ref{general_random}. 
In Section~\ref{corollaries} we show some of the consequences of Observation~\ref{edge_bound}.
Finally, in Section~\ref{questions} we list some unsolved questions which naturally came up during the research.

\section{Proof of Theorem~\ref{main_comp}} \label{main}

We start the proof with the cases in which the graph $H$ is equal to either $K_1 \oplus K_1$, 
$K_2$ or $K_2 \oplus K_1$.
The graph $G$ has $K_2$-avoiding coloring if and only if $G$ is the empty graph.

\smallskip

\noindent {\bfseries Proof of Lemma~\ref{K_2+K_1}.}
If $H=K_2 \oplus K_1$ then any two color classes in an $H$-avoiding coloring spans either 
a complete or empty bipartite graph. (In the special case if $H=K_1 \oplus K_1$ then there
can be only complete bipartite graphs between any two color classes.) Let us define a
binary relation $\rho$ such that for $x, y \in V(G)$ we have $x \rho y$ iff $N(x)=N(y)$.
Obviously $\rho$ is an equivalence relation, and the equivalence classes induced by $\rho$
are exactly the color classes of $G$ in the unique $K_2 \oplus K_1$-avoiding coloring of
$\chi_{K_2 \oplus K_1}(G)$ classes. \bizveg

\smallskip

Kr\'al, Kratochv\'\i l, Tuza and Woeginger \cite{KKTW} studied the hardness of 
coloring $H$-free graphs. They gave a complete description of the problem in the theorem follows:

\begin{theorem}[Kr\'al-Kratochv\'\i l-Tuza-Woeginger~\cite{KKTW}]  \label{2ktw} The problem 
$H$-{\rm F}{\scshape ree} {\rm C}{\scshape oloring} is polynomial-time solvable if $H$ is 
an induced subgraph of $P_4$ or of $P_3 \oplus K_1$, and NP-complete for any other $H$.
\end{theorem}  
   
Combining Theorem~\ref{2ktw} and Observation~\ref{triv}, one gets immediately that the
computation of $\chi_H(G)$ is NP-complete if $\chi(G)$ is also NP-complete for $H$-free $G$. On the other
hand, the polynomial-time computability of $\chi$ for $H$-free graphs does not imply the
same for $\chi_H$. Among the polynomial cases of Theorem~\ref{2ktw} we have checked 
already the graphs $K_1 \oplus K_1$, $K_2$ and $K_2 \oplus K_1$. Somehow against 
intuition, the computation of $\chi_H$ is NP-complete for the remaining $H=P_3$ and 
$H=P_4$ cases according to Lemmas~\ref{P_3} and \ref{P_4}.

\smallskip

\noindent {\bfseries Proof of Lemma~\ref{P_3}.}
We need to show $L_5=\{G: \chi_{P_3}(G) \leq 5\}$ is NP-complete language. We use reduction 
from the language $$L_{3, 2}=\{T: T \: \mbox{\rm is a 3-uniform hypergraph}, \chi(T) \leq 2\},$$ 
which is a well-known NP-complete problem. Let $T$ be an instance, that is $T \in L_{3, 2}$.

We need to assign a graph $G_T$ to $T$ such that $\chi_{P_3}(G_T) \leq 5$ if and only if
$\chi(T) \leq 2$. It turns out that the greatest difficulty is to associate the colorings
of the graph $G_T$ and the hypergraph $T$. The color of a vertex $t$ of $T$ cannot be encoded 
in one vertex $x_t$ of $G_T$, since the gadgets constructed in $G_T$ that enforce the good 
coloring of the edges of $T$ containing $t$ would interfere with each other. The solution
is to repeat the actual color of vertex $t$ at least as many times as the number of edges of 
$T$ that contain $t$.
For simplicity we repeat the color of any vertex $t$ a total of $m$ times, where $m$ is the
number of edges in $T$, and read out the color of $t$ at most once from each place.

The graph $G_T$ will consist of an $n \times m$ matrix of pentagons, in which the $i$-th row codes the
color of the $i$-th vertex in $T$. To assess the coloring of the $j$-th edge of $T$, the 
$j$-th column of this matrix is read. The usual types of gadgets are used in $G_T$ representing and
evaluating the edges of $T$, see Figure~\ref{fig1}.

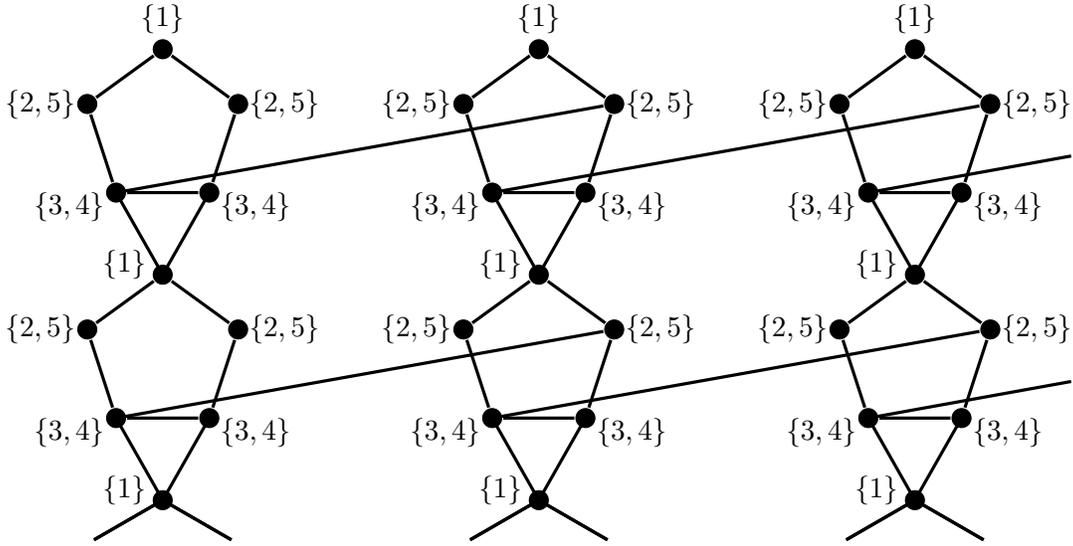
\begin{figure}
	\pgfdeclarelayer{background}
	\pgfdeclarelayer{foreground}
	\pgfsetlayers{background,main,foreground}
	\def\lgradius{30pt}
	\def\smradius{4pt}
									
	\def\xstr{5.0}
	\def\ystr{1.5}
	\begin{center}
		\begin{tikzpicture}[scale=1]
			\useasboundingbox (-1.5*\xstr,-2.75*\ystr) rectangle (1.5*\xstr,2.2*\ystr);
			
			\foreach \i in {-1,0,1}{
				\coordinate (\i;-3;0) at ($(\xstr*\i,\ystr*-3)+(90:\lgradius)$);
				\foreach \j in {-3,-1,1}{
					\foreach \k in {0,1,2,3,4}{
						\coordinate (\i;\j;\k) at ($(\xstr*\i,\ystr*\j)+(90+360*\k/5:\lgradius)$);
					}
				}
				\coordinate (\i;-3;1) at ($(\i*\xstr,-3*\ystr)+(30:\lgradius)$);
				\coordinate (\i;-3;4) at ($(\i*\xstr,-3*\ystr)+(150:\lgradius)$);
			}				
			\coordinate (2;-1;4) at ($(1*\xstr,-1*\ystr)+(-10:2*\lgradius)$);
			\coordinate (2;1;4) at ($(1*\xstr,1*\ystr)+(-10:2*\lgradius)$);
			
			\foreach \i in {-1,0,1}{
				\begin{pgfonlayer}{main}
					\draw [white, thin, fill=black] (\i;-3;0) circle(\smradius);
				\end{pgfonlayer}
				\begin{pgfonlayer}{foreground}
					\node[label={[xshift=2pt,yshift=4pt]180:{\small $\{1\}$}}] at (\i;-3;0) {};
					\node[label={[xshift=2pt,yshift=4pt]180:{\small $\{1\}$}}] at (\i;-1;0) {};
					\node[label={[label distance=-2pt]90:{\small $\{1\}$}}] at (\i;1;0) {};
				\end{pgfonlayer}
				\foreach \j in {-1,1}{
					\foreach \k in {0,1,2,3,4}{
						\begin{pgfonlayer}{main}
							\draw [white, thin, fill=black] (\i;\j;\k) circle(\smradius);
						\end{pgfonlayer}
					}
					\begin{pgfonlayer}{foreground}
						\node[label={[xshift=4pt,yshift=0pt]180:{\small $\{2,5\}$}}] at (\i;\j;1) {};
						\node[label={[xshift=4pt,yshift=6pt]205:{\small $\{3,4\}$}}] at (\i;\j;2) {};
						\node[label={[xshift=-4pt,yshift=6pt]335:{\small $\{3,4\}$}}] at (\i;\j;3) {};
						\node[label={[xshift=-4pt,yshift=0pt]0:{\small $\{2,5\}$}}] at (\i;\j;4) {};
					\end{pgfonlayer}
					\foreach \k in {0,1,2,3,4}{
						\begin{pgfonlayer}{background}
							\pgfmathtruncatemacro{\next}{mod(\k+1,5)};
							\draw [black,very thick] (\i;\j;\k) -- (\i;\j;\next);
						\end{pgfonlayer}
					}
					\pgfmathtruncatemacro{\below}{\j-2};
					\draw [black,very thick] (\i;\j;2) -- (\i;\below;0) -- (\i;\j;3);
					\pgfmathtruncatemacro{\next}{\i+1};
					\draw [black,very thick] (\i;\j;2) -- (\next;\j;4);
					\draw [black,very thick] (\i;-3;1) -- (\i;-3;0) -- (\i;-3;4);
				}
			}
		\end{tikzpicture}
	\end{center}
	\caption{The graph $G_T$ with the possible coloring, proof of Lemma~\ref{P_3}.}
	\label{fig1}
\end{figure}
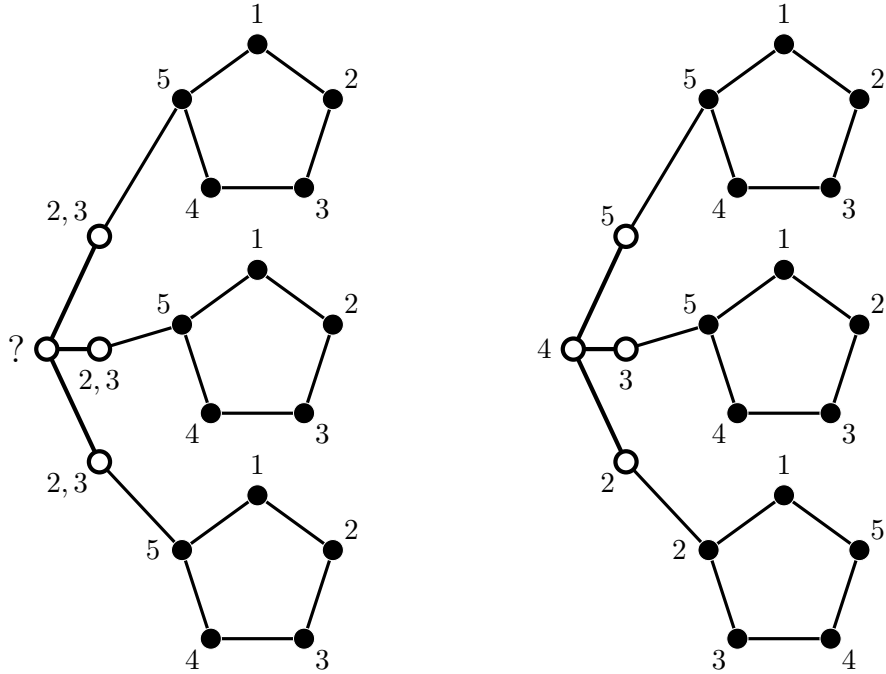
\begin{figure}
	\pgfdeclarelayer{background}
	\pgfdeclarelayer{foreground}
	\pgfsetlayers{background,main,foreground}
	
	\def\lgradius{30pt}
	\def\smradius{4pt}
									
	\def\xstr{7}
	\def\ystr{3}
	\begin{center}
		\begin{tikzpicture}[scale=1]
			\useasboundingbox (-0.5*\xstr,-1.5*\ystr) rectangle (1.2*\xstr,1.6*\ystr);
			
			\foreach \i in {0,1}{
				\foreach \j in {-1,0,1}{
					\foreach \k in {0,1,2,3,4}{
						\coordinate (\i;\j;\k) at ($(\xstr*\i,\ystr*\j)+(90+360*\k/5:\lgradius)$);
					}
				}
			}				
			\coordinate (A0) at (-0.4*\xstr,0.0*\ystr);
			\coordinate (A1) at (-0.3*\xstr,0.5*\ystr);
			\coordinate (A2) at (-0.3*\xstr,0.0*\ystr);
			\coordinate (A3) at (-0.3*\xstr,-0.5*\ystr);
			\coordinate (B0) at (0.6*\xstr,0.0*\ystr);
			\coordinate (B1) at (0.7*\xstr,0.5*\ystr);
			\coordinate (B2) at (0.7*\xstr,0.0*\ystr);
			\coordinate (B3) at (0.7*\xstr,-0.5*\ystr);
			
			\foreach \i in {0,1}{
				\foreach \j in {-1,0,1}{
						\begin{pgfonlayer}{foreground}
							\node[label={[yshift=0pt]90:{\small $1$}}] at (\i;\j;0) {};
						\end{pgfonlayer}
					\foreach \k in {0,1,2,3,4}{
						\begin{pgfonlayer}{main}
							\draw [white, thin, fill=black] (\i;\j;\k) circle(\smradius);
						\end{pgfonlayer}
					}
					\foreach \k in {0,1,2,3,4}{
						\begin{pgfonlayer}{background}
							\pgfmathtruncatemacro{\next}{mod(\k+1,5)};
							\draw [black,very thick] (\i;\j;\k) -- (\i;\j;\next);
						\end{pgfonlayer}
					}
				}
			}
			\foreach \i in {0,1}{
				\foreach \j in {0,1}{
					\begin{pgfonlayer}{foreground}
						\node[label={[xshift=4pt,yshift=-4pt]135:{\small $5$}}] at (\i;\j;1) {};
						\node[label={[xshift=4pt,yshift=4pt]225:{\small $4$}}] at (\i;\j;2) {};
						\node[label={[xshift=-4pt,yshift=4pt]315:{\small $3$}}] at (\i;\j;3) {};
						\node[label={[xshift=-4pt,yshift=-4pt]45:{\small $2$}}] at (\i;\j;4) {};
					\end{pgfonlayer}
				}
			}
			\begin{pgfonlayer}{foreground}
				\node[label={[xshift=0pt,yshift=0pt]180:{\small $5$}}] at (0;-1;1) {};
				\node[label={[xshift=4pt,yshift=4pt]225:{\small $4$}}] at (0;-1;2) {};
				\node[label={[xshift=-4pt,yshift=4pt]315:{\small $3$}}] at (0;-1;3) {};
				\node[label={[xshift=-4pt,yshift=-4pt]45:{\small $2$}}] at (0;-1;4) {};
			\end{pgfonlayer}
			\begin{pgfonlayer}{foreground}
				\node[label={[xshift=0pt,yshift=0pt]180:{\small $2$}}] at (1;-1;1) {};
				\node[label={[xshift=4pt,yshift=4pt]225:{\small $3$}}] at (1;-1;2) {};
				\node[label={[xshift=-4pt,yshift=4pt]315:{\small $4$}}] at (1;-1;3) {};
				\node[label={[xshift=-4pt,yshift=-4pt]45:{\small $5$}}] at (1;-1;4) {};
			\end{pgfonlayer}
			\draw [black, ultra thick, fill=white] (A0) circle(\smradius);
			\draw [black, ultra thick, fill=white] (B0) circle(\smradius);
			\foreach \x in {1,2,3}{
				\begin{pgfonlayer}{main}
					\draw [black, ultra thick, fill=white] (A\x) circle(\smradius);
					\draw [black, ultra thick, fill=white] (B\x) circle(\smradius);
				\end{pgfonlayer}
				\begin{pgfonlayer}{background}
					\draw [black,ultra thick] (A0) -- (A\x);
					\draw [black,ultra thick] (B0) -- (B\x);
				\end{pgfonlayer}
			}
			\begin{pgfonlayer}{background}
				\draw [black,very thick] (A1) -- (0;1;1);
				\draw [black,very thick] (B1) -- (1;1;1);
				\draw [black,very thick] (A2) -- (0;0;1);
				\draw [black,very thick] (B2) -- (1;0;1);
				\draw [black,very thick] (A3) -- (0;-1;1);
				\draw [black,very thick] (B3) -- (1;-1;1);
			\end{pgfonlayer}
			\begin{pgfonlayer}{foreground}
				\node[label={[xshift=4pt,yshift=-4pt]135:{\small $2,3$}}] at (A1) {};
				\node[label={[xshift=0pt,yshift=0pt]270:{\small $2,3$}}] at (A2) {};
				\node[label={[xshift=4pt,yshift=4pt]225:{\small $2,3$}}] at (A3) {};
				\node[label={[xshift=0pt,yshift=0pt]180:{\large\bfseries $?$}}] at (A0) {};
				\node[label={[xshift=4pt,yshift=-4pt]135:{\small $5$}}] at (B1) {};
				\node[label={[xshift=0pt,yshift=0pt]270:{\small $3$}}] at (B2) {};
				\node[label={[xshift=4pt,yshift=4pt]225:{\small $2$}}] at (B3) {};
				\node[label={[xshift=0pt,yshift=0pt]180:{\small $4$}}] at (B0) {};
			\end{pgfonlayer}
		\end{tikzpicture}
	\end{center}
	\caption{The two cases of colorings with gadgets, proof of Lemma~\ref{P_3}.}
	\label{fig2}
\end{figure}

Before examining the coloring of $G_T$, let us examine a $P_3$-avoiding good coloring of just $C_5$, since $C_5$ is the main 
building block of our construction. The vertices are referenced clockwise. If the first vertex 
is colored by 1, the second by 2, the third vertex color can be neither 2, because of adjacency,
or 1 since it would create a two colored $P_3$. So, without loss of generality, the first three
vertices are colored 1, 2, 3 respectively. The fourth vertex needs the fourth color. It cannot be colored by 2 or 3 as before.
If it would be colored by 1, the first, fifth and fourth vertices would form a two-colored 
$P_3$. Finally, the fifth vertex needs to be colored 5, since 1 and 4 are colors of adjacent
vertices, while 2 or 3 would create two-colored $P_3$'s.

Assuming that the $n$ vertices of $T$ are $x_1, \dots, x_n$ and the $m$ edges are $e_1, \dots, e_m$, 
the graph $G_T$ is first constructed by taking an $n \times m$ matrix with a $C_5$ in each position $(i,j)$. The $C_5$ 
in the $(i, j)$ position will be referred to as $C_{i, j}$. Second, connect the third and fourth vertices of 
$C_{i, j}$ to the first vertex of $C_{i+1, j}$ for $i=1, \dots, n-1$, $j=1, \dots n$, 
and similarly from $C_{n, j}$ to $C_{1, j+1}$ for $j=1, \dots, n-1$. Third, draw edges from the fourth vertex of $C_{i, j}$ to the second vertex of $C_{i, j+1}$ for $i=1, \dots, n$ and $j=1, \dots, m-1$. See Figure~\ref{fig1}. 

Without loss of generality, any $P_3$-avoiding five-coloring should use color 1 at the first vertex of any
$C_5$, should use the colors 2 and 5 in the second and fifth vertices (although in any order) and
the colors 3 and 4 in the third and fourth vertices (again their order is arbitrary). 

Furthermore, is easy to verify that a proper $P_3$-avoiding 
five-coloring must use the same order of colors 2 and 5 within a row, while the order of 
2 and 5 can be arbitrary for each row. We will use the $i$-th row to code the color of the vertex 
$x_i$ of the hypergraph $T$. However, when we read this ``value," each $C_5$ is read
only once.

Finally, the gadgets realizing the edges of $T$ are $m$ copies of $K_{1, 3}$. Let 
$e_\ell$ be $\{x_p, x_q, x_r \}$ and connect the leaves of the $\ell$-th $K_{1, 3}$ to the fifth
vertex of a yet unused $C_5$ in the $p$-th, $q$-th and $r$-th rows, respectively. The colors the 
vertices of $e_\ell$ receives are the color of the fifth vertices of $C_5$-s which were
connected to the leaves of the representing $K_{1, 3}$.  

Let us check if proper five-colorings of the construction and proper two-colorings of $T$
correspond to each other. If, for $e_\ell$, the vertices in the graph coloring all receive the color, say 
5, then the leaves of the representative $K_{1, 3}$ can be colored 2 or 3. One of these colors
appears two times, and it results in a two-colored $P_3$ in the graph coloring. If $e_\ell$
is colored properly, say 5, 5, 2, then the connected vertices in the representative $K_{1, 3}$
may get the colors 2, 3, 5. Giving color 4 to the 3-degree vertex of the representative 
$K_{1, 3}$ we get a proper $P_3$-avoiding five-coloring of $G$. See Figure~\ref{fig2}.

\smallskip

\noindent {\bfseries The case $ \chi_{P_3}(P_k) \leq 3$.}
If $G$ has a vertex of degree at least three, then $P_3(G)$ 
has a clique of size at least four, and by Observation~\ref{P_3close}, 
$\chi_{P_3}(G) \geq 4$. If all vertices have degree at most two, then the components 
of $G$ are paths and cycles. The components can be colored independently of each other 
in that case, so $G$ has a $P_3$-avoiding 3-coloring if and only if all components 
have. For all $k \in \mathbb N$, $\chi_{P_3}(P_k) \leq 3$, we just repeat the 
pattern $1, 2, 3, 1, 2, 3 \dots $ starting from one of the ends. The same can be 
(and must be) done for $C_k$ by specifying a starting vertex. However, it 
is successful only if $k \equiv 0 \mod 3$. \bizveg   
\bigskip

\noindent {\bfseries Proof of Lemma~\ref{P_4}.}
As in the proof of Lemma~\ref{P_3}, we use a reduction from the language $L_{3, 2}$, the 
two-coloring of 3-uniform hypergraphs. Having an instance $T \in L_{3, 2}$ with vertex set 
$x_1,\ldots,x_n$, $n\geq 4$, the reduction to a $P_4$-avoiding 3-coloring of a graph $G_T$ 
goes as follows. 
To each vertex $x_i$ of $T$ we create a pair of vertices $x_i, x_i'$ and have the 
edge $(x_i, x_i')$. An additional special vertex $z$ is adjacent to each $x_i$ and to each $x_i'$. 

\begin{figure}
	\pgfdeclarelayer{background}
	\pgfdeclarelayer{foreground}
	\pgfsetlayers{background,main,foreground}
	
	\def\lgradius{30pt}
	\def\smradius{4pt}
									
	\def\xstr{2}
	\def\ystr{2}
	\begin{center}
		\begin{tikzpicture}[scale=1]
			\useasboundingbox (-5.5,-3.5) rectangle (5.5,3.5);
			
			\coordinate(X) at (-4,1);
			\coordinate(X') at (-2,1);
			\coordinate(Z) at (-3,3);
			\foreach \j in {1,2,3}{
				\coordinate(x;1;\j) at (1,{4-\j});
				\coordinate(x;2;\j) at (-5,-\j);
				\coordinate(x;3;\j) at (1,-\j);
				\coordinate(a;1;\j) at (2,{4-\j});
				\coordinate(a;2;\j) at (-4,-\j);
				\coordinate(a;3;\j) at (2,-\j);
				\coordinate(b;1;\j) at (3,{4-\j});
				\coordinate(b;2;\j) at (-3,-\j);
				\coordinate(b;3;\j) at (3,-\j);
				\coordinate(c;1;\j) at (4,{4-\j});
				\coordinate(c;2;\j) at (-2,-\j);
				\coordinate(c;3;\j) at (4,-\j);
			}
			\coordinate(w1;1) at (4.8,2.5);
			\coordinate(w1;2) at (-1.2,-1.5);
			\coordinate(w1;3) at (4.8,-1.5);
			\coordinate(w2;1) at (4.8,1.5);
			\coordinate(w2;2) at (-1.2,-2.5);
			\coordinate(w2;3) at (4.8,-2.5);
						
			\begin{pgfonlayer}{main}
				\draw [white, thin, fill=black] (X) circle(\smradius);
				\draw [white, thin, fill=black] (X') circle(\smradius);
				\draw [white, thin, fill=black] (Z) circle(\smradius);
			\end{pgfonlayer}
			\foreach \i in {1,2,3}{
				\foreach \j in {1,2,3}{
					\begin{pgfonlayer}{main}
						\draw [black, ultra thick, fill=white] (x;\i;\j) circle(\smradius);
						\draw [white, thin, fill=black] (a;\i;\j) circle(\smradius);
						\draw [white, thin, fill=black] (b;\i;\j) circle(\smradius);
						\draw [white, thin, fill=black] (c;\i;\j) circle(\smradius);
					\end{pgfonlayer}
				}
				\begin{pgfonlayer}{main}
					\draw [white, thin, fill=black] (w1;\i) circle(\smradius);
					\draw [white, thin, fill=black] (w2;\i) circle(\smradius);
				\end{pgfonlayer}
			}
			\begin{pgfonlayer}{foreground}
				\node[label={[yshift=-2pt]270:{\small $x_i$}}] at (X) {};
				\node[label={[yshift=0pt]270:{\small $x_i'$}}] at (X') {};
				\node[label={[yshift=0pt]90:{\small $z$}}] at (Z) {};
				\node[label={[xshift=0pt]180:{\small $x_p$}}] at (x;1;1) {};
				\node[label={[xshift=6pt,yshift=-3pt]135:{\small $a_1$}}] at (a;1;1) {};
				\node[label={[xshift=6pt,yshift=-3pt]135:{\small $b_1$}}] at (b;1;1) {};
				\node[label={[xshift=6pt,yshift=-3pt]135:{\small $c_1$}}] at (c;1;1) {};
				\node[label={[yshift=2pt]270:{\small $1$}}] at (c;1;1) {};
				\node[label={[xshift=0pt]180:{\small $x_q$}}] at (x;1;2) {};
				\node[label={[xshift=6pt,yshift=-3pt]135:{\small $a_2$}}] at (a;1;2) {};
				\node[label={[xshift=6pt,yshift=-3pt]135:{\small $b_2$}}] at (b;1;2) {};
				\node[label={[xshift=6pt,yshift=-3pt]135:{\small $c_2$}}] at (c;1;2) {};
				\node[label={[yshift=2pt]270:{\small $2$}}] at (c;1;2) {};
				\node[label={[xshift=0pt]180:{\small $x_r$}}] at (x;1;3) {};
				\node[label={[yshift=2pt]270:{\small $2$}}] at (x;1;3) {};
				\node[label={[xshift=6pt,yshift=-3pt]135:{\small $a_3$}}] at (a;1;3) {};
				\node[label={[yshift=2pt]270:{\small $1$}}] at (a;1;3) {};
				\node[label={[xshift=6pt,yshift=-3pt]135:{\small $b_3$}}] at (b;1;3) {};
				\node[label={[yshift=2pt]270:{\small $3$}}] at (b;1;3) {};
				\node[label={[xshift=6pt,yshift=-3pt]135:{\small $c_3$}}] at (c;1;3) {};
				\node[label={[yshift=2pt]270:{\small $1$}}] at (c;1;3) {};
				\node[label={[xshift=0pt]0:{\small $w_1$}}] at (w1;1) {};
				\node[label={[xshift=0pt]0:{\small $w_2$}}] at (w2;1) {};
				\node[label={[yshift=2pt]270:{\small $2$}}] at (x;2;1) {};
				\node[label={[yshift=2pt]270:{\small $1$}}] at (a;2;1) {};
				\node[label={[yshift=2pt]270:{\small $3$}}] at (b;2;1) {};
				\node[label={[yshift=2pt]270:{\small $2$}}] at (c;2;1) {};
				\node[label={[yshift=2pt]270:{\small $1$}}] at (x;2;2) {};
				\node[label={[yshift=2pt]270:{\small $2$}}] at (a;2;2) {};
				\node[label={[yshift=2pt]270:{\small $3$}}] at (b;2;2) {};
				\node[label={[yshift=2pt]270:{\small $1$}}] at (c;2;2) {};
				\node[label={[yshift=2pt]270:{\small $1$}}] at (x;2;3) {};
				\node[label={[yshift=2pt]270:{\small $2$}}] at (a;2;3) {};
				\node[label={[yshift=2pt]270:{\small $3$}}] at (b;2;3) {};
				\node[label={[yshift=2pt]270:{\small $1$}}] at (c;2;3) {};
				\node[label={[xshift=0pt]0:{\small $3$}}] at (w1;2) {};
				\node[label={[xshift=0pt]0:{\small $2$}}] at (w2;2) {};
				\node[label={[yshift=2pt]270:{\small $1$}}] at (x;3;1) {};
				\node[label={[yshift=2pt]270:{\small $2$}}] at (a;3;1) {};
				\node[label={[yshift=2pt]270:{\small $3$}}] at (b;3;1) {};
				\node[label={[yshift=2pt]270:{\small $1$}}] at (c;3;1) {};
				\node[label={[yshift=2pt]270:{\small $2$}}] at (x;3;2) {};
				\node[label={[yshift=2pt]270:{\small $1$}}] at (a;3;2) {};
				\node[label={[yshift=2pt]270:{\small $3$}}] at (b;3;2) {};
				\node[label={[yshift=2pt]270:{\small $2$}}] at (c;3;2) {};
				\node[label={[yshift=2pt]270:{\small $1$}}] at (x;3;3) {};
				\node[label={[yshift=2pt]270:{\small $2$}}] at (a;3;3) {};
				\node[label={[yshift=2pt]270:{\small $3$}}] at (b;3;3) {};
				\node[label={[yshift=2pt]270:{\small $1$}}] at (c;3;3) {};
				\node[label={[xshift=0pt]0:{\small $3$}}] at (w1;3) {};
				\node[label={[xshift=0pt]0:{\small $3$}}] at (w2;3) {};
			\end{pgfonlayer}
			\draw[very thick] (Z) -- (X) -- (X') -- (Z);
			\foreach \i in {1,2,3}{
				\foreach \j in {1,2,3}{
					\begin{pgfonlayer}{background}
						\draw[ultra thick] (x;\i;\j) -- (a;\i;\j);
						\draw[very thick] (a;\i;\j) -- (b;\i;\j) -- (c;\i;\j);
					\end{pgfonlayer}
				}
			}
			\foreach \i in {1,2,3}{
				\begin{pgfonlayer}{background}
					\draw[very thick] (c;\i;1) -- (w1;\i) -- (c;\i;2) -- (w2;\i) -- (c;\i;3);
				\end{pgfonlayer}
			}
		\end{tikzpicture}
	\end{center}
	\caption{The two cases of colorings with gadgets, proof of Lemma~\ref{P_4}.}
	\label{fig3}
\end{figure}
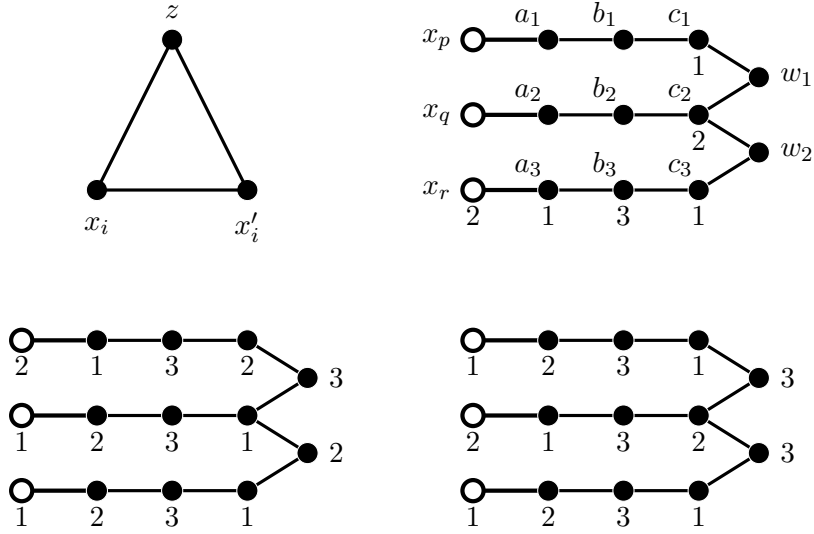

For each hyperedge $e_\ell=\{x_p, x_q, x_r\}$ in $T$, we define a gadget as follows. 
Take three disjoint copies of $P_3$, $a_i, b_i, c_i$ for $i=1, 2, 3$, and vertices $w_1, w_2$, 
and draw the edges $(c_1, w_1), (c_2, w_1), (c_2, w_2)$ and $(c_3, w_2)$. Finally we
set the gadget by drawing the edges $(a_1, x_p)$, $(a_2, x_q)$ and $(a_3, x_r)$. 

We claim that a $G$ has a $P_4$-avoiding 3-coloring if and only if $T$ has a $P_4$-avoiding 
2-coloring. We may assume vertex $z$ is colored by 3, so $x_i$s are colored 1 or 2, both
in the coloring of $G$ and $T$. 

\smallskip

If the vertices of an edge $e_\ell =\{x_p, x_q, x_r\}$
all receive the same color, say 2, then in the gadget associated to $e_\ell$ the
vertices $a_1, a_2, a_3$ must receive the color 1. (Indeed, if say $a_1$ would be colored
by 3, then take an $x_i \not \in \{x_p, x_q, x_r\}$. Either $x_i$ or $x_i'$ has the color 
2, say it is $x_i$. But $a_1, x_p, z, x_i$ is a 2-colored $P_4$.) The vertices $b_1, b_2$ 
and $b_3$ must get color 3, since if, say $b_1$, is of color 2, then $b_1, a_1, x_p, x_p'$
would induce a 2-colored $P_4$. If $c_2$ has color 1, then both $w_1$ and $w_2$ have 
color 2, since otherwise $w_1$ or $w_2$, $c_2$, $b_2$, $a_2$ would be a 2-colored $P_4$.
But in that case, the color of $c_3$ could be only 1, inducing 2-colored $P_4$ on 
the vertices $c_3, w_2, c_2, w_1$.

If $c_2$ has color 2, and at least one of $c_1$ or $c_3$ has color 1, assume $c_1$,
then $w_1$ must be colored 3. But then $w_2, c_1, b_1, a_1$ would be a 2-colored $P_4$.
Finally, if all $c_1, c_2, c_3$ has color 2, then $w_1$ and $w_2$ must have different
colors in order to avoid the 2-colored $P_4$ on $c_1, w_2, c_2, w_2$. But if, say
$w_1$ has color 3, then we see a 2-colored $P_4$ on the vertices $c_1, w_1, c_2, b_2$.

\smallskip

For the other direction, assume that $e_\ell=\{x_p, x_q, x_r\}$ received two colors 
in the hypergraph coloring. Without loss of generality, we may assume two vertices are 
colored 1, 
and one with 2. The vertex colored by 2 is either on one of the side, $x_p, x_r$ or 
the in the middle, $x_q$. Let us say $x_p$ has color 2 and $x_q, x_r$ received color 1. 
Then the coloring extends to the gadget of $e_\ell$ by coloring $a_1, c_2, c_3$ by 1,
$c_1, a_2, a_3, w_2$ by 2, and $b_1, b_2, b_3, w_1$ by 3. 
If $x_q$ has color 2 and $x_p, x_r$ have color 1, then the extension is 
$c_1, a_2, c_3$ is of color 1, $a_1, c_2, a_3$ is of color 2, and 
$b_1, b_2, b_3, w_1, w_2$ is of color 3. Notice that in both cases all $b$ types 
vertices received color 3, which ``insulates" the gadgets from each other, so the
defined 3-coloring is a $P_4$-avoiding one. \bizveg

\section{Proof of Theorem~\ref{2K_2poly3}} \label{polysec}

\noindent {\bfseries Proof of the case $H=2K_2$.}
To see if $\chi_{2K_2}(G) \leq 3$ for a given graph $G$, first we check if $G$ contains
$4K_2$ as an induced subgraph. This requires no more than $O(n^4)$ time. If $G$ does contain a $4K_2$,
then $\chi_{2K_2}(G) \geq 4$, since between two color classes there can be only one
of those four independent edges. Assume $G$ does not contain $4K_4$, and recall a 
result of Farber, Hujter and Tuza \cite{FHT}:

\begin{theorem} [Farber-Hujter-Tuza~\cite{FHT}] \label{poly} If the graph $G$ does not 
contain $(t+1)K_2$ as an induced subgraph, then the number of maximal independent 
sets in $G$ is at most $\binom{n/t}{t}$.
\end{theorem}

The following ideas are well-known and perhaps motivated Theorem~\ref{poly}.
The set $\mathcal M$ of all maximal independent sets can be found by, for example, 
a DFS tree algorithm, and can be listed in no more than $O(n^2|{\mathcal M}|)$ time. 
The decision problem of whether $\chi(G) \leq k$ can be solved by checking if there is $k$-set from $\mathcal M$ covering the vertex set of $G$. This still can be done in 
$O(\binom{|\mathcal M|}{k})$ time.

Applying Theorem~\ref{poly} to $G$, $|\mathcal M| \leq \binom{n/3}{3} < n^3/162$, so for a possible 
3-coloring we have to check a configuration of size no larger than $O(n^9)$. A configuration
consists of three maximal independent sets $X_1, X_2$ and $X_3$. First, $\cup_i X_i$
should contain all vertices of $G$. If this holds, it readily gives a 3-coloring, 
however it is not necessarily $2K_2$-avoiding. Indeed we are looking for 
$Y_i \subset X_i$ for $i=1, 2, 3$ such that $\cup_i Y_i$ contains all vertices of
$G$, $Y_i \cap Y_j = \emptyset$ if $i \not = j$, and the partition $\{Y_1, Y_2, Y_3\}$ 
is $2K_2$-avoiding.

We can assume that $\cap_i X_i=\emptyset$, if not, these vertices are isolated, and 
can be assigned to any $Y_i$ in the end. Then we start with the sets 
$Y_1:=X_1 \setminus (X_2 \cup X_3)$, $Y_2:=X_2 \setminus (X_1 \cup X_3)$ and
$Y_3:=X_3 \setminus (X_1 \cup X_2)$ and try to put the leftover vertices into those.  
The triple $\{Y_1, Y_2, Y_3\}$ should be $2K_2$-avoiding, otherwise we discard
the configuration. Then we have to decide, for example, if a vertex $x \in X_1 \cap X_2$ 
should be put in $Y_1$ or $Y_2$. If either placement would give a $2K_2$ with the 
set $Y_3$, we discard the configuration; if only one, we place it to the other; 
if none, we decide about it later. 

At the end of this process we have disjoint sets 
$Y_1, Y_2, Y_3$ that are $2K_2$-avoiding, $Y_i \subset X_i$, and the vertices of
$R_{1,2}:=(X_1 \cap X_2) \setminus (Y_1 \cup Y_2)$ can be placed both $Y_1$ or $Y_2$ 
(and same for $R_{1,3}$ and $R_{2,3}$). Let us construct a conflict graph on $R_{1, 2}$ and for
other indices do similarly. For $x, y \in R_{1, 2}$ there is an edge $(x, y) \in E(R_{1, 2})$ 
if $x$ and $y$ cannot be placed to $Y_1$. (That is, they induce a $2K_2$ to $Y_3$. It 
means $x$ and $y$ could not be placed in $Y_2$ either.) It is easy to see that if 
all those conflict graphs $R_{i, j}$, $i \not = j$ are bipartite, then all vertices
can be placed and we are ready. Otherwise the configuration is to be discarded and
we have to move to the next one. If none of the configurations can be formed to be a
$2K_2$-avoiding 3-coloring, then $\chi_{2K_2} > 3$. \bizveg

\section{Proof of Theorem \ref{general_random}} \label{random}

\noindent{\bfseries Proof of Proposition~\ref{upper_bound}.} 

Let $G$ have a coloring with part sizes $s_1,s_2,\ldots,s_\chi$ and $s_1$ the largest. 
First, further partition each color class arbitrarily into subparts of size at most $k_1-1$.
The number of parts is 
\begin{align*}
	\sum_{i=1}^{\chi} \left\lceil\frac{s_i}{k_1-1}\right\rceil 	\leq 	\sum_{i=1}^{\chi} \left(\frac{s_i}{k_1-1} + \frac{k_1-2}{k_1-1}\right) 
											= 	\frac{n}{k_1-1} + \frac{k_1-2}{k_1-1}\,\chi ,
\end{align*}
which is an upper bound that holds regardless of the value of $k_2$. 

Second, if $k_2\geq 3$, partition each color class except the largest arbitrarily into subparts of size at most $k_2-1$. 
The number of parts is 
\begin{align*}
	1 + \sum_{i=2}^{\chi} \left\lceil\frac{s_i}{k_2-1}\right\rceil 	&\leq 	1 + \sum_{i=2}^{\chi} \left(\frac{s_i}{k_2-1} + \frac{k_2-2}{k_2-1}\right) 
												= 		1 + \frac{n-s_1}{k_2-1} + \frac{k_2-2}{k_2-1}(\chi-1) \\
												&\leq 	1 + \frac{n-n/\chi}{k_2-1} + \frac{k_2-2}{k_2-1}(\chi-1) .
\end{align*}

Third, if $k_2=2$, color $G$ by giving the largest independent set one color and every other vertex an individual color. 
The number of parts is $n-\alpha+1$. 
Trivially, each of these partitions is an $H$-free coloring. 
All three combined bounds give the result in the proposition. 

\bizveg

\noindent{\bfseries Proof of Theorem~\ref{general_random}.} 

To obtain the upper bound, in the case of $k_1\geq 3$, we use Proposition~\ref{upper_bound} together with the result from Bollob\'as~\cite{Bol} that, \whp~$\chi\left(G(n,p)\right)=(1+o(1))\frac{n}{2\log_d n}$. Hence,
\begin{align*}
	\chi_H(G) 	\leq 	\frac{n}{k_1-1} + O\left(\frac{n}{\log_d n}\right) .
\end{align*}

In the case of $k_1=2$, the upper bound comes from Proposition~\ref{upper_bound} together with the result from Bollob\'as and Erd\H{o}s~\cite{BE} that, \whp~$\alpha\left(G(n,p)\right) = 2\log_d n-2\log_d\log n+O(1)$. Hence, 
\begin{align*}
	\chi_H(G) 	\leq 	n - 2\log_d n + O\left(\log_d\log n\right) .
\end{align*}

\smallskip

Now we proceed to the lower bound. An {\em $(\ell;k)$-complex} is a family of $\ell$ disjoint independent sets, each of size $k$, $A_1,\ldots,A_\ell$ such that each pair $(A_i,A_j)$, $1\leq i<j\leq\ell$ induces a graph that has no induced copy of $H$. 
The key to the proof is to show that for certain values of $k$ and $\ell=\ell(n)$, the probability that a $(\ell;k)$-complex exists goes to zero. 

If no $(\ell;k)$-complex exists, then whenever there is a coloring with color classes of size $n_1,\ldots,n_t\geq k$ it is the case that $\sum_{i=1}^t \lfloor n_i/k\rfloor < \ell$. Thus,
\begin{align*}
	\frac{1}{k}\sum_{i=1}^t n_i - \frac{k-1}{k} t 	\leq \sum_{i=1}^t \left\lfloor\frac{n_i}{k}\right\rfloor 	&< 	\ell \\
	\sum_{i=1}^t n_i 															&< 	k\ell + (k-1) t,
\end{align*}
while the leftover vertices are in color classes of size at most $k-1$.  So, if there are $t$ color classes of size at least $k$, then
\begin{align}
	\chi_{H}(G) 	\geq 	t + \frac{n - k \ell - (k-1) t}{k-1} = 	\frac{n}{k-1} - \frac{k}{k-1}\ell . \label{eq:general_random}
\end{align}

For the graph $H$, let $Q=Q(H,p)$ be the probability that a $k_1\times k_1$ random bipartite graph has no induced copy of $H$. 
Taking the product over all $\binom{\ell}{2}$ pairs $(A_i,A_j)$ and multiplying by the probability that each $G[A_i]$ induces an independent set, we obtain:
\begin{align}
	\Pr\left[\exists\text{ an }(\ell;k_1)\text{-complex}\right] 	
	&= 	\frac{(n)_{k_1\ell}}{\ell! (k_1!)^\ell} Q^{\binom{\ell}{2}} (1-p)^{\ell\binom{k_1}{2}} \label{eq:pr_bound} \\
	&< 	\left[\left(\sqrt{\frac{e}{\ell}}\right) n Q^{(\ell-1)/(2k_1)} (1-p)^{(k_1-1)/2}\right]^{k_1\ell} , \nonumber
\end{align}
which is obtained from the inequalities $(n)_{k_1\ell}\leq n^{k_1\ell}$, $\ell!\geq (\ell/e)^{\ell}$, and $k_1!\geq 1$. 

Let $C'=C'(H,p)=\frac{2k_1}{\log (1/Q)}$. For $n$ sufficiently large, if $\ell > C' \log n$, then the probability in \eqref{eq:pr_bound} goes to zero.
By \eqref{eq:general_random},  it is the case that \whp
\begin{align*}
	\chi_{H}(G) 	\geq 	\frac{n}{k_1-1} - C' \frac{k_1}{k_1-1} \log n . 
\end{align*}
Thus the general lower bound is satisfied. 

In the special case where $H=2K_2$, we observe that $Q(2K_2,p)=1-2p^2(1-p)^2$. 
\begin{align*}
	\Pr\left[\exists\text{ an }(\ell;2)\text{-complex}\right] 	
	&= 	\frac{(n)_{2\ell}}{\ell! 2^\ell} Q^{\binom{\ell}{2}} (1-p)^{\ell} \\
	&< 	\left[\left(\sqrt{\frac{e}{2\ell}}\right) n Q^{(\ell-1)/4} (1-p)^{1/2}\right]^{2\ell} .
\end{align*}
If $\ell>\frac{4\log n}{\log (1/Q)} - \frac{4\log\log n}{\log (1/Q)} + \log\log\log n$ and $n$ is sufficiently large, then \whp~no $(\ell;2)$-complex exists.
By \eqref{eq:general_random},  it is the case that \whp
\begin{align*}
	\chi_{2K_2}(G) 	\geq 	n - (1-o(1)) \frac{8 \log n}{\log (1/Q)} , 
\end{align*}
where $Q=1-2p^2(1-p)^2$.

\bizveg

Similar results can be obtained as long as $\min\{p,1-p\}=\omega\left(\frac{\log n}{n}\right)$ but express our results in the case where $p$ is a fixed constant.  \\

\section{Examples for Observation~\ref{edge_bound}} \label{corollaries}

An easy consequence of Observation~\ref{edge_bound} is as follows:

\begin{corollary} \label{P_n1}
	\begin{align*}
		\chi_{2K_2}(P_n) 	\geq 	\sqrt{2 \left\lceil\frac{n-1}{3}\right\rceil + \frac{1}{4}} + \frac{1}{2} .
	\end{align*}
\end{corollary}

A more refined argument gives the value of $\chi_{2K_2}(P_n)$ as follows:

\begin{corollary} \label{P_n2}
	If $k$ is the least integer that satisfies 
	\begin{align*}
		\left\lfloor\frac{k+1}{2}\right\rfloor (k-2) 	\geq 	\left\lceil\frac{n-1}{3}\right\rceil	
	\end{align*}	
	then $\chi_{2K_2}(P_n) = k$.
\end{corollary}

\noindent{\bfseries Proof.} 
Let $\ell = \left\lceil\frac{n-1}{3}\right\rceil$. In particular, this means $n-1\leq 3\ell\leq n+1$.

Choose $k$ to be the least integer so that
\begin{align}
	\left\{
		\begin{array}{lll}
			\binom{k}{2} 			&\geq 	\ell+1, 	& \mbox{ if $k$ is odd;} \\
			\binom{k}{2}-\frac{k}{2}+1 	&\geq 	\ell+1, 	& \mbox{ if $k$ is even.}
		\end{array}
	\right.
\label{eq:kl}
\end{align}

This value is chosen because the longest Eulerian trail in $K_k$ has $\binom{k}{2}$ edges if $k$ is odd and has $\binom{k}{2}-\frac{k}{2}+1$ edges if $k$ is even. The latter case occurs when a matching of size $k/2-1$ is removed from $K_k$.

Let $a_1,a_2,\ldots,a_{\ell},a_{\ell+1}$ be an Eulerian trail in $K_k$. 
Enumerate the vertices of $P_n$ as $1,2,\ldots,n$. Let the coloring of the vertices of $P_n$ be $f : [n] \rightarrow \{a_1,a_2,\ldots,a_{\ell},a_{\ell+1}\}$, defined as follows:
	\begin{align*}
		f(1) = f(3) 					&= 	a_1 . \\
		f(3i-4) = f(3i-2) = f(3i) 		&= 	a_{i}, 		&&\mbox{ for $i=2,\ldots,\ell-1$.} \\
		f(3\ell-4) = f(3\ell-2) 			&= 	a_{\ell} . \\
		f(3\ell) 					&= 	a_{\ell}, 		&&\mbox{ if $3\ell\leq n$.} \\
		f(3\ell-1) 					&= 	a_{\ell+1} . \\
		f(3\ell+1)					&= 	a_{\ell+1}, 	&&\mbox{ if $3\ell+1\leq n$.}					
	\end{align*}
Note that $3\ell+2\geq n+1$ by our choice of $\ell$, so there are no other vertices to color.

For all $j\in\{1,\ldots,\ell-1\}$ each pair of color classes $a_ja_{j+1}$ induces a $P_4$ plus some isolated vertices. The pair $a_{\ell}a_{\ell+1}$ also induces a path plus isolated vertices and the path is either $P_2$, $P_3$, or $P_4$, depending on the remainder of $n$ modulo 3. Therefore, this coloring is a $2K_2$-free coloring.

\smallskip

To see that equality holds, consider a $2K_2$-free $k$-coloring of $P_n$. Every pair of color classes either induces an empty graph or a graph whose only nontrivial component is $P_2$, $P_3$, or $P_4$. Thus, the edges of $P_n$ can be partitioned according to which unique pair of colors induce a particular subpath. Furthermore, subpaths that share a vertex must share a color. Thus, we can construct an auxiliary graph $\Gamma$ on $\{1,\ldots,k\}$ where $ij$ is an edge if and only if the pair of colors $\{i,j\}$ induces a path on at least two vertices. Because consecutive small paths must share a vertex and hence a color, the edges of $\Gamma$ form a trail on $K_k$. Since each edge of $\Gamma$ corresponds to at most 3 edges of $P_n$, the number of edges in $P_n$ is at most three times the length of a longest trail in $K_k$. That is,
\begin{align*}
	n-1 	\leq 		3\cdot \left\{ 	\begin{array}{ll}
								\binom{k}{2}, 				& \mbox{ if $k$ is odd;} \\
								\binom{k}{2}-\frac{k}{2}+1, 	& \mbox{ if $k$ is even.}
							\end{array} \right.
\end{align*}

\smallskip

Finally, we return to \eqref{eq:kl} and observe that the condition on $k$ is equivalent to 
\begin{align*}
	\left\lfloor\frac{k+1}{2}\right\rfloor (k-2) 	\geq 	\ell .
\end{align*}
The statement then follows. See Table~\ref{Pn_table} for small values of $\chi_{2K_2}(P_n)$.
\bizveg

\begin{table}
	\begin{tabular}{|r||c|c|c|c|c|c|c|c|c|} \hline
		$n\in$ 		& [2,4] 	& [5,7] 	& [8,13] 	& [14,28] 	& [29,37] 	& [38,61] 	& [62,73] 	& [74,106] 	& [107,121] \\ \hline
		$\chi_{2K_2}(P_n)$ 	& 2		& 3 		& 4 		& 5 		& 6 		& 7 		& 8 		& 9 			& 10 \\ \hline
	\end{tabular}
	\caption{$\chi_{2K_2}(P_n)$ for small values of $n$.}
	\label{Pn_table}
\end{table}

\begin{corollary} \label{matching}
\begin{align*}
		\chi_{2K_2}\left(\frac{n}{2}\cdot K_2\right) 	= \left\lceil \sqrt{n+\frac{1}{4}}+\frac{1}{2} \right\rceil .
	\end{align*}
\end{corollary}

\noindent{\bfseries Proof.}
If $k = \chi_{2K_2}\left(\frac{n}{2}\cdot K_2\right)$, then 
$\binom{k}{2} \geq \frac{n}{2}$ by Observation~\ref{edge_bound}, which establishes a lower bound of 
$\left\lceil \sqrt{n+\frac{1}{4}}+\frac{1}{2} \right\rceil$. 
This is also sufficient in that the vertices of each edge of $\frac{n}{2}\cdot K_2$ can be 
independently assigned to distinct endvertices of $K_k$ where $k$ is the value given. \bizveg

Bounds on some other graphs are given as follows:

\begin{corollary} \label{tree}
Let $n$ be odd and let $T$ be the tree formed when each edge of $K_{1,(n-1)/2}$ is subdivided (by a vertex) exactly once.
	\begin{align*}
		\chi_{2K_2}(T) 	= \left\lceil \sqrt{n-\frac{3}{4}}+\frac{1}{2} \right\rceil .
	\end{align*}
\end{corollary}

The proof for Corollary~\ref{tree} is nearly identical to that of Corollary~\ref{matching} so it is left to the reader.

\begin{corollary} \label{cube}
Let $Q_d$ be the $d$-dimensional hypercube. Then
	\begin{align*}
		\chi_{2K_2}(Q_2) 	&= 		2, \\
		\chi_{2K_2}(Q_3) 	&= 		4, \\
	\end{align*}
and if $d \geq 4$ then
	\begin{align*}
		\chi_{2K_2}(Q_d) 	&\geq 	\sqrt{\frac{d}{2 d - 1} 2^d} + \frac{1}{2}  \\
						&= 		\sqrt{\frac{n}{2} \frac{1}{1 - 1/(2\lg n)}} + \frac{1}{2} .
	\end{align*}
\end{corollary}

\noindent{\bfseries Proof.}
The case of $Q_2$ is trivial. It takes some small work to show that for $Q_3$ there is no $Q_3$-avoiding 
$3$-coloring and there is a $Q_3$-avoiding $4$-coloring. In the case where $d \geq 4$, we use Observation~\ref{edge_bound}, 
where $k = \chi_{2K_2}(Q_d)$, 
	\begin{align*}
		(2 d - 1) \binom{k}{2} 	\geq 	d 2^{d-1} ,
	\end{align*}
because the only graphs that can occur between pairs of color classes are double-stars 
(at most $2d-1\geq 7$ edges) and $K_{2,2}$ ($4$ edges). 

\bizveg

\begin{definition} \label{projective_plane}
Let $p$ be a prime power and let $G(p)$ be the bipartite graph on $n = 2 (p^2 + p + 1)$ 
vertices defined by the projective plane of order $p+1$. That is, there are $p^2 + p + 1$ 
points and $p^2 + p + 1$ lines and a point is adjacent to a line if and only the point is 
in the line in the projective plane. This graph is $(p + 1)$-regular with no $K_{2,2}$. 
\end{definition}

\begin{corollary} \label{chrproj} If $G(p)$ is the graph in Definition~\ref{projective_plane}
then,
	\begin{align*}
		\chi_{2K_2}(G(p)) 	&\geq 	\sqrt{\frac{2 (p^2 + p + 1) (p + 1)}{2 p + 1}} + \frac{1}{2} \\
					& \geq 		\sqrt{\frac{n}{2} + \frac{\sqrt{n}}{4 }} + \frac{1}{2} .
	\end{align*}
\end{corollary}

\noindent{\bfseries Proof.}
With $k = \chi_{2K_2}(G)$, we use the inequality
	\begin{align*}
		(2 d + 1) \binom{k}{2} \geq (p^2 + p + 1) (p + 1) .
	\end{align*}

\bizveg

\section{Questions} \label{questions}

The proof method of Theorem~\ref{2K_2poly3} suggests the following problem:

\begin{quest}\label{2K_2poly_k}
Is it polynomial time decidable if $\chi_{2K_2}(G) \leq k$ for any fixed $k \in \mathbb{N}$?
\end{quest}

The structure of $P_4$-free graphs is very nice, which allows easy computation of 
chromatic number, clique number etc. In fact $P_4$-free bipartite graphs, called {\em difference graphs}~\cite{HPS}, are well-studied. It is quite surprising that the function 
$\chi_{P_4}$ is NP-hard. 

\begin{quest}
Is it true that determining whether $\chi_{P_4}(G)\leq 3$ is NP-hard?
\end{quest}

\begin{quest}
	What is $\chi_{2K_2}(Q_d)$?
\end{quest}

\begin{quest}
	What is $\max\{\chi_{2K_2}(T) : v(T)=n,  T \textrm{ is a tree}\}$?
\end{quest}

\begin{quest}
	Is it true that $\chi_{2K_2}(G)=n-(1+o(1))2\log_2n$ \whp~if $G\sim G(n,1/2)$? 
\end{quest}



\end{document}